\documentclass[11pt,leqno]{amsart}
\usepackage{latexsym}
\usepackage{amssymb,amsbsy}
\usepackage[mathscr]{euscript}
\usepackage{amsmath}
\usepackage{amsthm}
\usepackage{amscd}
\usepackage{amsfonts}
\newtheorem{theorem}{Theorem}[section]
\newtheorem{lemma}[theorem]{Lemma}
\newtheorem{corollary}[theorem]{Corollary}

\theoremstyle{definition}

\newtheorem{question}{Question}
\theoremstyle{remark}
\newtheorem{remark}[theorem]{Remark}

\newcommand{\h}{{\mathbb H}}

\newcommand{\C}{{\mathbb C}}

\newcommand{\dsp}{\displaystyle}

\newcommand{\R}{{\mathbb R}}
\newcommand{\T}{{\mathbb T\,}}
\newcommand{\I}{{\mathbb I}}

\title{On unitary representability of topological groups}
\author{Jorge Galindo}
\address{Departmento de Matem\'aticas \\
         Universidad Jaume I \\
        Campus Riu Sec, 12071\\
        Castell\'on
        Spain}
\email{jgalindo@mat.uji.es} \subjclass[2000]{43A35,46B99, 22A10,
54H11, 54E35.}

\keywords{Unitary group, positive definite, free abelian topological
group, free locally convex space, free Banach space, unitarily
representable, uniform embedding, Schwartz space}

\date{}

\thanks{Research  partly supported by the Spanish Ministry of Science (including
FEDER funds), grant MTM2004-07665-C02-01. The foundations of this
paper were laid during  the authors stay at the University of Ottawa
supported by a Generalitat Valenciana grant CTESPP/2004/086.}
\begin{document}
 \maketitle
\begin{abstract}
 We prove that the  additive group $(E^\ast,\tau_k(E))$ of an
$\mathscr{L}_\infty$-Banach space $E$, with  the topology
$\tau_k(E)$ of uniform convergence on compact subsets of $E$,  is
 topologically
isomorphic to a subgroup of the unitary group of some Hilbert space
(is \emph{unitarily representable}). This is the same as proving
that the topological group $(E^\ast,\tau_k(E))$ is uniformly
homeomorphic to a subset of $\ell_2^\kappa$ for some $\kappa$.

As an immediate consequence,  preduals of commutative von Neumann
algebras or duals of commutative $C^\ast$-algebras are unitarily
representable in the topology of uniform convergence on compact
subsets. The unitary representability of  free locally convex spaces
(and thus of free Abelian topological groups) on compact spaces,
follows as well.

The above facts cannot be extended to noncommutative von Neumann
algebras or general Schwartz spaces.
 \end{abstract}
\section{Introduction}\label{1}
We address in this paper the problem of determining which
topological groups embed in the unitary group of a Hilbert space.
This problem is relevant in the program of extending Harmonic
Analysis beyond locally compact groups (see \cite{pest00} or
    \cite{pest05}, see also  \cite{shte94}), and strongly
 related to other much studied problems like the uniform
embeddability  of metric or Banach spaces in $\ell_2$, see
\cite[Chapter 7]{benylind}, or the universality of orbit equivalence
relations induced by actions of the unitary group, see
\cite{gao,gaopest}.

 When a topological group can be embedded in the unitary
group of a Hilbert space (unitary groups will be assumed to carry
the strong operator topology) we say that it is \emph{unitarily
representable}. Positive definite functions play a prominent
r\^{o}le in the unitary representability of a topological group.
This is  mainly because the  (strong operator) topology of  the
unitary group is determined
 by positive definite functions. As a matter
of fact, a topological group is
 unitarily representable precisely when its topology is generated  by
its continuous positive definite functions as happens,  for
instance, with  discrete groups (characteristic functions of points
are linear combinations of positive definite functions), see Lemma
\ref{char2}. All locally compact groups are actually   unitarily
representable, their regular representation $\lambda$ establishes a
topological isomorphism $\lambda\colon G\to \mathcal{U}(L_2(G))$. If
we go beyond the class of locally compact groups, regular
representations are no longer at reach. Megrelishvili \cite{megr00}
has proven
 that isometry groups of separable Hilbert spaces
  embed uniformly in $\ell_2$ and, therefore, that for
  separable Banach spaces (or even
Polish groups), unitary representability implies uniform
embeddability in $\ell_2$.   This shows how strongly   the problem
of unitary representability is linked with another problem that has
been the object of deep research by many authors: the uniform
classification of
 Banach
spaces.  Both problems are  indeed equivalent for Abelian groups,
see Theorem \ref{fchar}. For non-Abelian groups it remains unknown
whether these two problems remain  equivalent.

 We will consider here
the unitary representability of additive groups of Banach spaces
equipped  with topologies of uniform convergence on compact sets,
the natural topology in the duality theory of topological groups.
Denote by $\tau_k(E)$ the topology of uniform convergence on compact
subsets of a topological space $E$. We prove in Section \ref{3} that
for an arbitrary Banach $\mathscr{L}_\infty$-space $E$,
 the
additive group of  $(E^\ast,\tau_k(E))$ embeds isomorphically in a
power of $\ell^1$ and is therefore  unitarily representable. The
results stated in the second paragraph of the abstract then follow.
Since $(E^\ast,\tau_k(E))$ is always a Schwartz locally convex space
(i.e. a locally convex space that can be obtained as projective
limit of Banach spaces wit compact linking maps) and so are other
classes  of unitarily representable groups, such as nuclear spaces
it is natural to study whether Schwartz spaces must be unitarily
representable. We analyze this question in Section 4 and focus on
Schwartz topologies on the Banach space $c_0$. We find that
pre-images of  neighbourhoods of 1 in $\C$ under continuous positive
definite function on $c_0$ are  always quite large. As a
consequence, we obtain that positive-definite functions do not
define either the norm topology or the topology of uniform
convergence on compact subsets of $\ell_1$. Since the latter is a
Schwartz topology,   we conclude that Schwartz spaces need not be
unitarily representable.  This example also shows that dual groups
of  Banach spaces may fail to be unitarily representable.

We finish this introduction with  some words on notation and
terminology.

 Unless specifically stated all topological groups will
be Abelian. Given a finite measure $\mu$, the topological vector
space of all complex-valued measurable functions will be denoted by
$L_0(\mu)$, $L_0(\mu,\T)$ will stand for the subset of $L_0(\mu)$
consisting of $\T$-valued functions. Under multiplication,
$L_0(\mu,\T)$ becomes a topological group. We will always assume
that $L_0(\mu)$ carries the topology of convergence in measure.
Sometimes we will find useful to specify the space $X$ on which
$\mu$ is defined and use $L_0(X,\mu)$ instead of $L_0(\mu)$.

By an uniform embedding between topological groups $T\colon G \to H$
we mean an injective mapping such that both $T$ and $T^{-1}$ are
uniformly continuous.

 By $\tau_k(X)$ we will denote the topology of uniform convergence on the compact subsets of a
topological space $X$. If $E$ is a Banach space, the term \emph{dual
group of $E$}, in symbols $\widehat{E}$,   will refer to the
additive group of $(E^\ast,\tau_k(E))$.  In the duality theory of
topological groups, the symbol $\widehat{E}$ is usually reserved for
the character group of $E$, i. e. the group of continuous
homomorphisms into the unit circle $\T$ with the topology of uniform
convergence on the compact subsets of $E$. It is easy to see that
this group is topologically isomorphic to $(E^\ast,\tau_k(E))$.

The  ball of radius $\delta$ centered in the identity  of a
metrizable group $G$  will be denoted by $B_{G,\delta}$, in  case
$\delta=1$, we will simply write $B_G$.
\section*{Acknowledgements}
Concerning this paper, I owe my thanks to several people that have
contributed  to it in  different senses.

The origin of this paper is in the author's discussions with V.
Pestov at the University of Ottawa in  the Spring of 2004. The
author is indebted to V. Pestov for proposing  him the present
problem and for providing the first arguments in favour of the
unitary representability of $A([0,1])$.

Some  observations of  V. Tarieladze were very useful. M.
Megrelishvili has contributed by carefully reading several drafts of
this paper and pointing out some inaccuracies. W. B. Johnson
suggested the final form of Theorem \ref{compactoper} and therefore
the final tenor of the paper.

 I am finally  indebted to V.  Uspenski\u\i for finding  a wrong
statement in a previous version of this paper and for sharing with
me an elegant, direct proof of the unitary representability of
$L(X)$ for completely regular $X$.
\section{Unitary representability of topological groups}\label{2}
 We  collect here some facts
available in the literature on unitary representability  of
topological groups.

Positive-definite functions provide an internal criterion for the
unitary representability of  a topological group, see Lemma
\ref{char2}, and therefore constitute an essential tool. A
complex-valued function $\phi$ defined on a topological group $G$ is
positive definite provided that for every finite subset $\{x_j
\}_{j=1}^n$ of $G$ and every
 collection of complex numbers $\{\lambda_{j}\}_{j=1}^n$,
\[ \sum_{i,j=1}^n \lambda_{i}\bar{\lambda_j}\phi(x_i\,x_j^{-1})\geq 0 .\]
The function $e^{-\| \cdot \|}$ defined on a Hilbert space, is a
positive definite function of special importance. Its positive
definiteness is usually proven as a consequence of Schoenberg's
theorem on exponentiation of so-called negative-definite functions,
see \cite[Chapter 3]{bergchriress}. This function suffices to note
one direction of the following well-known characterization of
 uniform embeddability  through positive definite functions.
\begin{lemma}
\label{char2} A topological group  $G$ is  unitarily representable
if and only if for each neighbourhood $U$ of the identity 0 of $G$,
there is a real-valued positive definite function $\phi_U$ with
$\phi_{_{U}}(0)=1$  such that
\[ \{ g \in G \colon |\phi_{_{U}}(g)-1|<1/2\}\subseteq U.\]
\end{lemma}
\begin{proof}
Sufficiency of the condition can be proved as in Theorem 2.1 of
\cite{gao}.

Assume now that $G$ is unitarily representable and let  $T \colon G
\to \mathcal{U}(\h)$ denote a topological isomorphism of $G$ into
the unitary group of some Hilbert space $\h$. Consider
$\widetilde{U}$ a neighbourhood of the identity in $\mathcal{U}(\h)$
such that $T(U)=\widetilde{U}\cap T(G)$. Since the topology of
$\mathcal{U}(\h)$ is the strong operator topology, there will be
some  $\varepsilon_i>0 $ and  $\xi_i \in \h$, $1\leq i\leq n$, such
that
\[\widetilde{U}\supset \bigcap_{i=1}^n  \left\{S\in \mathcal{U}(\h)\colon
\|S\xi_i-\xi_i\|\leq \varepsilon_i\right\}.\] Now define ${\dsp
\phi_{_{U}}= e^{-\log 2\sum_{i=1}^n \frac{1}{\varepsilon_i}
\|T(\cdot)\xi_i-\xi_i\|}}$.
 \end{proof}
 Schoenberg  \cite{scho}  not only realized  that
$e^{-\|\cdot\|}$ is a positive definite function on any $L_2(\mu)$,
he also discovered that raising the norm to a power $\alpha<1$ did
not change the positive definiteness of the norm and completely
classified $L_p$-spaces in terms of  the positive definiteness of
$e^{-\|\cdot\|}$.

This classification is contained in the following theorem that
summarizes several  known results about unitary representability of
Banach spaces.
\begin{theorem}[Aharoni, Maurey and Mityagin \cite{aharmaurmity},
Megrelishvili \cite{megr00}, Schoenberg \cite{scho}, see Chapter 8
of \cite{benylind}]\label{schol1} If a Banach space is unitarily
representable, then it can be isomorphically embedded in $L_0(\mu)$
for some measure $\mu$ and must have cotype 2. The Banach space
$L_p(\mu)$ is unitarily representable if and only if $1\leq p \leq
2$.
\end{theorem}
\begin{proof}
(Sketch) Megrelishvili \cite{megr00} proves that $U(\ell_2)$ embeds
uniformly in $\ell_2$ and Aharoni, Maurey and Mityagin prove in
\cite{aharmaurmity} that  a Banach space admits a uniform  embedding
in $\ell_2$ if and only if  it admits an \emph{isomorphic} embedding
in $L_0(\mu)$.  Since (1) $e^{-\|\cdot\|}$ is positive definite on
$L_p$ for $1\leq p \leq 2$ (Schoenberg), (2) every Banach space
isomorphic to a subspace of $L_0(\mu)$ has cotype 2 (\cite[Corollary
8.17]{benylind}) and (3) $L_p(\mu)$-spaces with $p>2$ do not have
cotype 2, all the assertions in the theorem follow.
\end{proof}
Gathering  the information of these two Sections, we obtain    the
following characterization of unitary representability.
\begin{theorem} \label{fchar}
Let $G$ be an Abelian topological group. The following assertions
are equivalent:
    \begin{enumerate}
\item $G$ is unitarily representable.
\item $G$ is topologically isomorphic to a subgroup of
 $L_0(X,\mu,\T)$, for some
compact $X$ and some  Borel measure $\mu$  on $X$.
 \item There is a uniform embedding of $G$
into   a power  $(\ell_1)^\kappa$ of the Banach space $\ell_1$.
\item There is a uniform embedding of $G$
into   a power $(\ell_2)^\kappa$ of the Banach space $\ell_2$.
\end{enumerate}
\end{theorem}
\begin{proof}
The equivalence between (1) and (2) is well-known and appears for
instance in \cite{gaopest} for Polish groups. For non-Polish groups
the argument is still valid.

Megrelishvili \cite{megr00} proves that the unitray group
$\mathcal{U}(\ell_2)$ embeds uniformly in $\ell_2$. Since the
unitary group $\mathcal{U}(\oplus_\kappa \h)$ of a
$\kappa$-dimensional Hilbert  space is topologically isomorphic to
the product $\prod_\kappa \mathcal{U}(\h)$,  this shows that (1)
implies (4).

That (1) follows from (4)
 is
 a pretty straightforward consequence of
 the results of \cite{aharmaurmity}, and  might have been
  known to some specialists. We  sketch the proof of this
  fact.

Assume  that (4) holds. Let $\ell_{2,i}$  be a copy  of $\ell_2$ for
every $i<\kappa$ and consider a  uniform embedding $T\colon G \to
\prod_{i<\kappa} \ell_{2,i}$. Let also $\{A_j: j<\kappa\}$ be an
enumeration of all finite subsets of $\kappa$ and define the
corresponding projections $\pi_j \colon\prod_{i<\kappa}
\ell_{2,i}\to\prod_{i\in A_j} \ell_{2,i}$. The norm of $\prod_{i\in
A_j} \ell_{2,i}$ ($\cong \ell^2$) will be denoted by $\|\cdot\|_j$.

We now   adapt  Corollary 3.6 of \cite{aharmaurmity} (see also
Chapter 8 of \cite{benylind}).

 Let  $M$ denote an invariant mean
on $G$ (i.e. a continuous functional on $\ell^\infty(G)$ with
$M(1)=1$ and $M(L_xf)=M(f)$, if $L_xf(y)=f(xy)$). Each  projection
$\pi_j$ defines   a continuous  positive definite function $\phi_j$
given by
\[ \phi_j(g) =M(f_g),\quad \text{ with } \quad
f_g(y)=e^{-\|\pi_jT(gy)-\pi_jT(y)\|_j^2}.\]To avoid confusions, we
have adopted here multiplicative notation for the group operation on
$G$ and additive notation for $\ell^2$. That $\phi_j$ is positive
definite follows from the invariance of $M$ and that $\phi_j$ is
continuous  follows from the uniform continuity of $\pi_j\,T$, see
Corollary 3.6 of \cite{aharmaurmity}. Since $T^{-1}$ is uniformly
continuous as well, every neighbourhood $U$ of the identity in $G$
determines a projection $\pi_j$ and an $\varepsilon>0$ such that
\[\{ g \colon \|\pi_j T(g)-\pi_j(\bar{a})\|_j<\varepsilon \}\subset
\{ g \colon g\in T^{-1}(\bar{a})\,U\} \quad \text{ for every }
\bar{a}\in T(G)\subset (\ell_2)^\kappa.\] Arguing again as in
Corollary 3.6 of \cite{aharmaurmity} it follows that there is
$\delta >0$  such that
\[ \{g \in G \colon |\phi_i(g)-1|<\delta\}\subset U.\]
Hence, $G$ is unitarily representable by Lemma \ref{char2}.

That assertions  (4) and (3) are  equivalent is indeed a  well-known
fact. That (3) implies (4) is proved in \cite{aharmaurmity}, since
$\ell_2$ embeds even isometrically in $\ell_1$, (4) implies (3).
\end{proof}
\begin{corollary}\label{unif}
An Abelian Polish  group is unitarily representable if and only if
it is uniformly embeddable in $\ell_2$.
\end{corollary}
\begin{remark}
 As M. Megrelishvili    \cite{megrxx} has indicated to the author,
Corollary \ref{unif} is true for Polish amenable groups, this
readily  follows from the proof of  the implication $(4)\implies
(1)$ of Theorem \ref{fchar}. According to this same letter,
Megrelishvili announced  that  (4) implies (1) in 2002. This was
however never published.
\end{remark}
\begin{corollary}\label{lh}
The additive group of a  Banach space is unitarily representable if
and only if it is both linearly isomorphic to a subspace of
$L_0(\mu)$ and topologically isomorphic to a subgroup of $L_0(\mu,
\T)$.
\end{corollary}

\section{Dual groups of $\mathscr{L}_\infty$-spaces}\label{3}
  We show in this section that the dual space  $E^\ast$  of
  an ${\mathscr L}_\infty$-space $E$ with the topology of uniform convergence
 on compact subsets of $E$, $\tau_k(E)$,
is always unitarily representable.

Recall  that  ${\mathscr L}_\infty$-spaces are Banach spaces whose
finite dimensional subspaces are close to $\ell_\infty^m$-spaces,
see for instance \cite[Appendix F]{benylind} for  the actual
definition. Here it suffices to say that $C(K)$-spaces are
$\mathscr{L}_\infty$ spaces for every compact space $K$.
\begin{lemma}[Theorem 4.3 and Remark 4.5 of
\cite{fonfetal}] Let $E$  and $Y$ be
     ${\mathscr L}_\infty$-spaces, $Y$ separable. For each
      compact subset $K$ of
  $E$ there is   a
  one-to-one compact operator $T\colon Y \to E$ with
  $K\subset T(B_Y)$.
  \end{lemma}
\begin{corollary}\label{compactoper}
Let $E$ be a $\mathscr{L}_\infty$ space and let $K\subset E$ be
compact. There is then a compact one-to-one operator $T_K\colon
c_0\to E$ with $K\subset T_K(B_{c_0})$.
\end{corollary}
 \begin{theorem} \label{fd}
The dual group  $\widehat{E}=(E^\ast,\tau_k(E))$ of a ${\mathscr
L}_\infty$-space $E$ embeds in a power of $\ell^1$, and hence is
unitarily representable.
\end{theorem}
\begin{proof}
Let $ \mathcal{K}(E)$ denote a set that is cofinal in the family of
all compact subsets of $E$ (ordered by inclusion). For each
$K\in\mathcal{K}(E)$ we consider  the compact operator $T_K\colon
c_{0,K}\to E$ defined on  a  copy $c_{0,K}$ of $c_0$ that is
provided by Corollary \ref{compactoper}.

Define now
\[ \Psi \colon E^\ast \to \prod_{K\in\mathcal{K}(E)}
 \ell_{1,K} \] as the product
$\Psi=\prod_{K\in\mathcal{K}(E)}  T_K^\ast$, where again
$\ell_{1,K}$ represents a copy of $\ell_1$. We now check that $\Psi$
is one-to-one, $\tau_k(E)$-continuous and open.

Suppose $\Psi(f)=0$, and take any $x\in E$. By considering some
$K\in \mathcal{K}(E)$ with $x\in K$, we see that $f(x)=0$, and
therefore that $\Psi$ is injective.

To see  that $\Psi$ is $\tau_k$-continuous, we must check that the
conjugate operators $T_K^\ast$ are all $\tau_k$(E)-continuous.
Taking into account that $T_K$ is a compact operator, we have that
$T_k(B_{c_{0,K}})\subset K_0$ for some $K_0\in \mathcal{K}(E)$.
Uniform convergence on $K_0$ will therefore imply uniform
convergence on $B_{c_{0,K}}$ (and thus norm-convergence on
$\ell_{1,K}$). Accordingly $T_K^\ast$ is $\tau_k$-continuous.

As to the openness of $\Psi$, it is sufficient  to observe that
uniform convergence on every $B_{c_{0,K}}$ implies
$\tau_k(E)$-convergence, since every $K\in \mathcal{K}(E)$ is
contained in the $T_K$-image of the  corresponding $B_{c_{0,K}}$.
\end{proof}
\begin{corollary}
If a topological group $G$ admits a  uniform embedding into
  $(L_1(X,\mu),\tau_k(C(X)))$
 for some compact $X$ and some  Borel measure $\mu$ on $X$, then $G$ is
 unitarily representable
 \end{corollary}

\begin{corollary}\label{vn}
Let $A$ denote a commutative Banach algebra, and let $\mathscr{M}$
denote a commutative von Neumann algebra with predual
$\mathscr{M}_\ast$. The additive groups of $A^\ast$ and
$\mathscr{M}_\ast$ are both unitarily representable for the
respective
 topologies of uniform convergence on compact subsets $\tau_k(A)$
 and $\tau_k(\mathscr{M})$.
\end{corollary}

 Theorem \ref{fd} can be directly applicable to free locally convex
spaces and free Abelian topological groups on compact spaces. If $X$
is a completely regular space, the free locally convex space $L(X)$
and the free Abelian topological group $A(X)$ on $X$ are obtained by
providing the free vector space (resp. the free Abelian group on
$X$) with a locally convex vector space topology (resp. a
topological group topology) such that every continuous function
$f\colon X \to E$ into a locally convex space (resp. a topological
grouop) can be extended to a continuous linear map $\bar{f}\colon
L(X)\to E$ (resp. to a continuous homomorphism), this is thus a
\emph{linearization} process. What is important for our purposes is
that   these free topological objects can  be realized as preduals
of $C(X,\R)$ and $C(X,\T)$
 respectively: $C(X,\R)$  appears as the
  space of continuous linear functionals on $L(X)$
  and  $C(X,\T)$ appears as the space of continuous
characters on $A(X)$. $L(X)$ carries  the topology of uniform
convergence on equicontinuous pointwise bounded subsets of,
$C(X,\R)$  and  $A(X)$  the topology of uniform convergence on
equicontinuous subsets of  $C(X,\T)$, see \cite{uspe83}, and
\cite{pest95}.

The dual group of  $C(X,\R)$ will be denoted as $M_c(X)$. Since the
topology of $L(X)$ is the topology of uniform convergence on
equicontinuous pointwise bounded (=relatively compact) sets, $L(X)$
is  a subgroup of $M_c(X)$ in a natural way. By
Tkachenko-Uspenski\u\i's theorem \cite{tkac83,uspe90}, $A(X)$ is a
topological subgroup of $L(X)$ and thus of $M_c(X)$. The following
is then an immediate corollary to Theorem \ref{fd}, it answers
Question 35 in \cite{pest06}, see also Question 6.10 in
\cite{megr06}
\begin{corollary} \label{compactm}
The additive group of  $M_c(K)$  and its subgroups $L(K)$ and $A(K)$
are unitarily representable for every compact Hausdorff space $K$.
\end{corollary}
\begin{remark}\label{all}
In response to a previous version of this paper, Uspenski\u\i
\cite{uspe06l} has provided a different argument
 for the unitary representabilty of $L(X)$,  valid for
\emph{any} completely regular space $X$. The preprint is now
available \cite{uspe06}.
\end{remark}
\begin{remark}
For a given metric space $(K,d)$, the free locally convex space over
$K$ can be obtained as the projective limit of the free Banach
spaces over $(K,d_\alpha)$ for some family of distances $d_\alpha$
compatible with $d$, see for instance \cite{gaopest} for the precise
meaning of the terms free Banach and free normed. For that reason,
it was somehow expected that Corollary \ref{compactm}
 should follow from an
 analogous statement for some
 free normed  spaces $FN(K,d_f)$.
  Notice however that Theorem \ref{fd} does
not allow such a deduction. Examples of free Banach spaces that are
not unitarily representable are actually easy to come by,
 $FN(\I^2)$, where $\I$ denotes the unit interval
  $\I=[0,1]$ is one such example.
  It is unclear what happens with  $FN(\I,d_f)$ for arbitrary
  metric transforms $d_f$,
  (i.e. metrics that arise after composing $d$ with an increasing.
  concave function $f\colon \I \to \I$).
 If $f_\alpha(t)=t^\alpha$, it is
known (\cite{cies60b} or \cite{weav}, see also the  recent complete
account of \cite{kalt04}) that $FB(\I,d_{f_\alpha})$ (the completion
of $FN(\I,d_f)$) is isomorphic to $\ell_1$ and hence unitarily
representable; the proof of this fact can be extended to other
metric  transforms but not to all of them, see \cite{cies60a}.
\end{remark}
\section{On the limits of the class of unitarily representable
groups} We try here to get  an idea of which are the limits of  the
class of unitarily representable groups. We have shown in Section 2
that dual groups of Banach $\mathscr{L}_\infty$-spaces are in that
class. Another distinguished family of unitarily representable
groups is the variety of nuclear groups that contains additive
groups of nuclear locally convex spaces and locally compact Abelian
groups, see  \cite{bana91}. The unitary representability of nuclear
groups was obtained in \cite{bana01}; in the particular case of
nuclear locally convex spaces, this follows directly from their
representation as projective limits of Hilbert spaces.

It is not clear how to extend  the class of unitarily representable
groups beyond nuclear groups or dual groups  of Banach
$\mathscr{L}_\infty$-spaces, as neighbouring families already
contain groups that are not unitarily representable.

One possible extension could involve a noncommutative version of
Corollary \ref{vn}. But Corollary \ref{vn} is based on embedding the
dual group $(A^\ast,\tau_k(A))$ of a commutative von Neumann algebra
$A$ in a product of preduals of commutative von Neumann algebras
(namely $\ell_1$'s). This point of view cannot be carried over to
the noncommutative case, as already preduals of noncommutative von
Neumann algebras may fail to be unitarily representable. The algebra
of trace-class operators $C_1$, which is the predual of the von
Neumann algebra $B(\ell_2)$ of all bounded operators on $\ell_2$, is
one example, see the Remark in page 194 of \cite{benylind}.

Another possible extension could involve Schwartz spaces. Dual
groups of $\mathscr{L}_\infty$-spaces and nuclear spaces both belong
to the class of Schwartz spaces (or to its topological
group-theoretic analog \emph{Schwartz groups} \cite{aussetal} for
the case of nuclear groups). Recall that a locally convex space is a
\emph{Schwartz space} if it can be represented as a projective limit
of Banach spaces with compact linking maps. Notwithstanding its
closeness to theses classes,  the class of Schwartz spaces, also
contains groups that are not unitarily representable. Our example
will be the additive group of $(c_0,\tau_k(\ell_1))$. This is in
fact a universal generator for the class of Schwartz spaces. Since
compact subsets of $\ell_1$ are contained in the closed convex hull
of null sequences, $\tau_k(\ell_1)$ can be replaced by the topology
of uniform convergence on null-sequences in $\ell^1$. We denote  by
$\mathscr{S}(c_0)$ the vector space $c_0$ equipped with this
topology.

Neighbourhoods of 0 in $\mathscr{S}(c_0)$, are determined by
sequences  of  numbers going to 0: given such a sequence
$\overline{\alpha}=(\alpha_n)_n$,  we consider the corresponding
neighbourhood
\[ U_{\overline{\alpha}}=\left\{ (x_n)_n\in c_0 \colon
|x_n|\leq \frac{1}{|\alpha_n|}, \text{ for every $n$}\right\}.\] The
sets $U_{\overline{\alpha}}$ (with $\alpha \in c_0$) constitute a
neighbourhood basis at 0 of $\mathscr{S}(c_0)$, see \cite{rand73}.

To see why $\mathscr{S}(c_0)$ is not unitarily representable we need
a couple of results on the structure of linear operators with values
in $L_0$.
\begin{theorem}[Nikishin factorization theorem, see for instance
Theorem 13 of \cite{dilw} or Proposition 8.16 of
\cite{aharmaurmity}] \label{niki} Let $X$ be a Banach space. Every
continuous linear operator  $T\colon X \to L_0$  factorizes through
$L_q$ for each  $0<q<1$, i.e. there are continuous linear operators
$S_1\colon X \to L_q$ and $S_2\colon L_q\to L_0$ such that
$S_2S_1=T$.
\end{theorem}
In the sequel we will make use of the standard unit vectors  $e_n
\in c_0$. The symbol  $e_n$
 will  denote,  as usual,  the sequence with one in the $n$th place
 and zero otherwise.
\begin{theorem}[Theorem 4.3 of \cite{kaltmont}]\label{orlicz}
For every bounded operator $T\colon c_0 \to L_p$, $0<p\leq 2$, $\sum
\|T(e_n)\|^2<\infty$.
 \end{theorem}
 \begin{theorem}\label{schw}
 The Schwartz space $\mathscr{S}(c_0)$ is not unitarily representable.
 \end{theorem}
  \begin{proof}
Let $\overline{\alpha}=(\alpha_n)_n$ denote the sequence with
$\alpha_n=1/\log(n+1)$ and consider the
$\mathscr{S}(c_0)$-neighbourhood of the identity
\[U_{\overline{\alpha}}=\left\{(a_n)_n\in c_0 \colon a_n\leq \log (n+1),\text{ for every $n$} \right\}.\]
We will  next see that there is no positive definite function $\phi$
on $c_0$ such that $\{(a_n)_n\in c_0 \colon
|\phi\bigl((a_n)_n\bigr)-1|<1/2\}$ is contained in
$U_{\overline{\alpha}}$. By Lemma \ref{char2}  this will prove the
theorem.

By Lemma  4.2 of \cite{aharmaurmity} each positive definite mapping
$\phi$  on a linear topological space $X$ induces a continuous
linear operator $T\colon X\to L_0$ in such a way that for every
$x\in X$, the characteristic function $\varphi_{_{T(x)}}$ of $T(x)$
satisfies $\varphi_{_{T(x)}}(t)=\phi(tx)$ for every $t\in(0,1)$.
Here  we have adopted the usual probabilistic terminology and by the
\emph{characteristic function} of   $f\in L_0(\mu)$ we mean the
complex-valued function $\varphi_f(t)=\int e^{i t f(x)}\, d\mu(x)$.
It is then easy to see that there is some $\delta>0$ such that
\begin{equation} T^{-1}(B_{L_0,\delta})\subset  \{(a_n)_n\in c_0 \colon
|\phi\bigl((a_n)_n\bigr)-1|<1/2\}.\label{cf} \end{equation}

Let $T\colon c_0 \to L_0$ and $\delta>0$ be  a continuous linear
operator and a real number respectively. In view of \eqref{cf}, it
will  suffice to see that $T^{-1}(B_{L_0,\delta})\nsubseteq
U_{\overline{\alpha}}$.

 Let $c_0\overset{S_1}{\to}{L_r}\overset{S_2}{\to}L_0$ be a factorization of the operator $T$ as
 in Theorem \ref{niki}. By Theorem \ref{orlicz},
$\sum_n \|S_1(e_n)\|_r^2=M<\infty$.

Let $\varepsilon >0$ be  such that $S_2(B_{L_r,\varepsilon})\subset
B_{{L_0,\delta}}$ and choose $k$ large enough so that
$\|S_1(e_k)\|_r<\frac{ \varepsilon}{ 2  \log(k+1)}$. Take  now  some
$(a_n)_n\in c_0$ with $a_n\leq \varepsilon/(2 M)$ if $n\neq k$ and
with $a_k=\frac{\varepsilon}{2 \|S_1(e_k)\|_r}$. Clearly
$(a_n)_n\notin U_{\overline{\alpha}}$ and
$\|S_1\bigl((a_n)_n\bigr)\|_r\leq \sum_n |a_n| \cdot
\|S_1(e_n)\|_r\leq \varepsilon$. We have thus that
$T\bigl((a_n)_n\bigr)=S_2S_1\bigl((a_n)_n\bigr)\in B_{L_0,\delta}$
while $(a_n)_n\notin U_{\overline{\alpha}}$.
\end{proof}
It immediately follows that dual groups of Banach spaces need not be
unitarily representable (and hence that  Theorem \ref{fchar} cannot
be thus generalized).
\begin{corollary}
The dual group of $\ell_1$, i.e. $(\ell_\infty,\tau_k(\ell_1)) $, is
not unitarily representable.
\end{corollary}
\begin{proof}
Just observe that $\mathscr{S}(c_0)$ embeds linearly as topological
subspace of $\widehat{\ell_1}=(\ell_\infty,\tau_k(\ell_1)) $.
\end{proof}
\begin{remark}
Theorem \ref{schw} shows that neighbourhoods of the identity for the
topology generated by continuous positive definite functions on
$c_0$ are quite large, and hence that this topology is rather weak.
This same reason also shows that $c_0$ does not embed uniformly in
$\ell_2$, a fact that was first proved by Enflo \cite{enflo}.
\end{remark}

 Most of the techniques of this paper are commutative in
nature,  it seems nevertheless worth to end this paper with a
mention to two natural questions that are left untouched here.
\begin{question}[Question 4.4 of \cite{megr00}]\label{q1}
If a  (non-Abelian) Polish topological group embeds uniformly in
$\ell_2$, must it be unitarily representable?\end{question}
\begin{question}[Question 35 of \cite{pest06}]
Is the free topological group on a compact space (or even $[0,1]$)
unitarily representable?
\end{question}
\def\cprime{$'$} \def\cprime{$'$}
  \def\polhk#1{\setbox0=\hbox{#1}{\ooalign{\hidewidth
  \lower1.5ex\hbox{`}\hidewidth\crcr\unhbox0}}}
  \def\polhk#1{\setbox0=\hbox{#1}{\ooalign{\hidewidth
  \lower1.5ex\hbox{`}\hidewidth\crcr\unhbox0}}}

\end{document}